\documentclass[11pt]{amsart} 
\usepackage{amsmath,amssymb,a4wide} 
\usepackage{color}
\newtheorem{theorem}{Theorem} 
\newtheorem{lemma}[theorem]{Lemma} 
\newtheorem{proposition}[theorem]{Proposition} 
\newtheorem{remark}[theorem]{Remark} 

\newtheorem{defi}[theorem]{Definition}

\newcommand{\pr}{\noindent  {\bf \textit{Proof : }}}

\newcommand{\cqfd}{{\nobreak\hfil\penalty50\hskip2em\hbox{}\nobreak\hfil
$\square$\qquad\parfillskip=0pt\finalhyphendemerits=0\par\medskip}}

\newcommand{\dis}{\displaystyle}

\newcommand{\tore}{\Omega}

\newcommand{\TD}{{\mathbb T}}

\let\eps\varepsilon 

\usepackage{color}

\title[About  the barotropic compressible quantum Navier-Stokes equations]
{About the barotropic compressible quantum Navier-Stokes equations}

\author[M. Gisclon]{M. Gisclon}
\address{Laboratoire de Math\'ematiques, CNRS UMR 5127,
Universit\'e de Savoie, 73376 Le Bour\-get-du-Lac, France; 
e-mail: gisclon@univ-savoie.fr}

\author[I. Violet]{I. Lacroix-Violet} 
\address{Laboratoire de Math\'ematiques, CNRS UMR 8524 
Universit\'e de Lille 1, Villeneuve d'Ascq, France;
e-mail: ingrid.violet@univ-lille1.fr}

\begin{document}


\begin{abstract}
In this paper we consider the barotropic compressible quantum Navier-Stokes equations with a linear density dependent viscosity and its limit when the scaled Planck constant vanish. Following recent works on degenerate compressible Navier-Stokes equations, we prove the global existence of weak solutions by the use of a singular pressure close to vacuum. With such singular pressure, we can use the standard definition of global weak solutions which also allows to justify the limit when the scaled Planck constant denoted by $\varepsilon$ tends to $0$. 
\end{abstract}


\maketitle
{\small\noindent 

}


\section{Introduction-Motivations}\label{sec.intro}

In this paper, the model of interest belongs to quantum fluid models. Such models can be used to describe superfluids \cite{LoMo93}, quantum semiconductors \cite{FeZhou93}, weakly interacting Bose gases \cite{Grant73} and quantum trajectories of Bohmian mechanics \cite{Wyatt05}. Recently some dissipative quantum fluid models have been derived. In \cite{GuJu04} the authors derived viscous quantum Euler models using a moment method in Wigner-Fokker-Planck equation. In \cite{BruMe09}, under some assumptions, using a Chapman-Enskog expansion in Wigner equation, the quantum Navier-Stokes equations are obtained. 

In this paper, we are interesting in the barotropic quantum Navier-Stokes equations  which read as, for $x \in \Omega$ and  $t>0$
\begin{equation} 
\label{eq1}
\left\{  
 \begin{array} {ll}
\partial_t n+\mbox{div}(n \, u) =  0, \\[2mm]
\displaystyle{\partial_t(n \, u)+\mbox{div}(n \, u \otimes u) +\nabla_x (p(n)) -2 \,  \varepsilon^2 \, n \nabla \left( \frac{\Delta \sqrt{n}}{\sqrt{n}} \right)  =   2 \, \nu \,  \mbox{div}(n \, D(u))},Ê\\[2mm]
n|_{t=0}=n_0, \,  \quad (n \, u)|_{t=0} = n_0 \, u_0,
\end{array}
\right.
\end{equation}
where the unknowns are the particle density $n$ and the particle velocity $u$. Here, $\Omega=\TD^d$ is the torus in dimension $d$ (in this article $1 \leq d \leq 3$), and $u \otimes u$ is the matrix with components $u_i \, u_j$. The function $p(n)=n^{\gamma}$ 
with $\gamma > 1$ is the pressure function and $D(u)$ stands for the symetric part of the velocity gradient, namely $D(u)=\displaystyle {(\nabla u +^t  \nabla u)}/{2}$. Finally, the physical parameters are the Planck constant $\varepsilon >0$ and the viscosity constant $\nu>0$. 

The quantum correction $(\Delta \sqrt{n})/\sqrt{n}$
 can be interpreted as a quantum potential,
the so-called Bohm potential, which is well known in quantum mechanics.  This  Bohm potential  arises from the fluid dynamical
formulation of the single-state Schr\"odinger equation. 
The non-locality of
quantum mechanics is approximated by the fact that the equations of state do not only depend
on the particle density but also on its gradient. 
These equations were employed to model field
emissions from metals and steady-state tunneling in metalÐinsulatorÐmetal structures and to simulate ultra-small semiconductor devices .

Formally, setting $\varepsilon=0$ in \eqref{eq1} leads to the compressible Navier-Stokes equation with density-dependent viscosity. Our goal is to perform the limit $\varepsilon$ tends to zero. To this end we need some global existence result on \eqref{eq1}.

Existence results for the stationary isothermal model in one space dimension were shown in \cite{JM12}. The
main mathematical difficulty, besides of the highly nonlinear structure of the
third-order quantum terms, is the proof of positivity (or non-negativity) of the
particle density.

In \cite{Jue10}, A. J\"ungel, proves the global existence of this system when the scaled Planck constant   is bigger than the viscosity constant ($\varepsilon > \nu$).
In \cite{Dong10}, J. Dong extends this result where the  viscosity constant is equal to  the scaled Planck constant ($\varepsilon = \nu$), and, in \cite{Jiang11}, F. Jiang shows that the result still holds when the viscosity constant is bigger than the scaled Planck constant ($\nu> \varepsilon$). Therefore, thanks to \cite{Jue10}, \cite{Dong10} and \cite{Jiang11},
we have the global existence for all physically interesting cases of the scaled Planck and viscosity constants. Note that the definition of global weak solutions used in \cite{Jue10}, \cite{Dong10} and \cite{Jiang11} follows the idea introduced in \cite{BrDeLi2003} by testing the momentum equation by $n\,\phi$ with $\phi$ a test function. Here the problem of such formulation is that it requires $\gamma>3$ for $d=3$ which is not a suitable assumption for physical case, and the estimates on the solution are badly dependent on $\varepsilon$ due to extra terms coming from the diffusive term tested against $n\phi$. 

Then, the goal of this paper is to prove existence of global weak solutions without the assumption $\gamma >3$ if $d=3$ and with uniform estimates 
allowing to perform the limit of \eqref{eq1} when $\eps$ tends to 0. To this end we will use an another formulation.  In more details, our proof of existence will be split in two steps. In a first time we construct some approximate solutions and in a second time, using some a priori estimates, we prove the stability of solutions. The construction of approximate solutions will follow exactly the ideas introduced by A.~J\"uengel in \cite{Jue10}. We will recall this approximate system which gives an other way to construct approximate solutions than the one introduced by \cite{MuPoZa2}  after  some hints given in \cite{BrDe2007}. The new lines in our paper,  compared to \cite{Jue10}, will be the asymptotic limit with respect to the regularized parameter. Indeed we replace the concept of global weak solutions obtained by multiplying the momentum equation by $n\phi$ by a more standard formulation which required the addition of an extra cold pressure as introduced in \cite{BrDe2007} and fully developed for chemical reactive flows in  \cite{MuPoZa2} and \cite{Za2013}. The idea is to prove global existence of weak solutions of the following system, for $x \in \Omega$ and $t>0$:
\begin{equation} 
\label{eq2}
\left\{  
 \begin{array}{ll}
\partial_t n+\mbox{div}(n \, u) = 0, \\[2mm]
\partial_t(n \, u)+\mbox{div}(n \, u \otimes u) +\nabla_x (p(n) + p_c(n))
        -2 \,  \varepsilon^2 \, n \nabla \left( \dis  \frac{\Delta \sqrt{n}}{\sqrt{n}} \right)  =  2 \, \nu \,  \mbox{div}(n \, D(u)),Ê\\[2mm]
n|_{t=0} =n_0, \,  \quad (n \, u)|_{t=0} =n_0 \, u_0 \, 
\end{array}
\right.
\end{equation}
where  $p_c$ is a singular continuous function namely a suitable increasing function satisfying 
$$\lim_{n \to 0} p_c(n) = + \infty$$
and called cold pressure. More precisely, we assume
\begin{equation} 
p'_c(n)=
\left\{  
 \begin{array} {l}
 c \, n ^{-4k-1}, \mbox{ for } n \leq 1, \, k>1, \\
 n^{\gamma-1}, \, \mbox{ for } n >1, \, \gamma >1
 \end{array}
\right.
\end{equation}
for some constant $c>0$.

\begin{remark}
As mentioned in \cite{BrDe2007}, the physical relevance of the compressible Navier-Stokes equations is very questionnable in regions where the density are close to vacuum: the medium is not only unlikely to be in a liquid or gas state (elasticity and plasticity have to be considered for such solid materials for which by the way low densities may lead to negative pressures), but also the rarefied regime of vanishing densities violates the assumptions on the mean free path of particles suitable for fluid models. The singular pressure is a way to modelize such plastification in fluid models. Of course, it could be interesting from a mathematical view-point  to understand if such singular pressure is necessary.
\end{remark}

\begin{remark}
Such degeneracy will help (at the level of stability) to conclude about the strong convergence of term of type$\sqrt {n} \, u$ and thus pass to the limit in the nonlinearity $n u \otimes u$. 
\end{remark}
We write the quantum term in a different form to pass to the limit: 
\begin{eqnarray} \nonumber
< n \nabla \left( \dis  \frac{\Delta \sqrt{n}}{\sqrt{n}} \right)  , \phi>
 & = &  < \nabla \left(\sqrt{n}  \Delta \sqrt{n}  \right)  , \phi>-<\frac{1}{\sqrt{n}} \nabla n \Delta \sqrt{n} , \phi> \\
  \nonumber
 & = & - < \sqrt{n} \Delta \sqrt{n}, \mbox{div} \phi>-2< \nabla \sqrt{n} \Delta \sqrt{n}, \phi> \\
 \nonumber
 & = &-< \nabla\left( \sqrt{n} \nabla \sqrt{n}\right),\mbox{div} \phi>+
< | \nabla \sqrt{n}|^2, \mbox{div} \phi> \\
\nonumber 
& & -2 <\mbox{div} (\nabla \sqrt{n} \otimes \nabla \sqrt{n}),\phi>+2 < \left(\nabla \sqrt{n} \cdot \nabla\right) \nabla \sqrt{n}, \phi> \\
\nonumber
 & = & < \sqrt{n} \nabla \sqrt{n}, \nabla \mbox{div} \phi> + < | \nabla \sqrt{n}|^2,\mbox{div} \phi>  \\
 \nonumber
  & & +2 < \nabla \sqrt{n} \otimes \nabla \sqrt{n}, \nabla \phi>
 + < \nabla | \nabla \sqrt{n}|^2,\phi> \\
 \nonumber
 & = &  < \sqrt{n} \nabla \sqrt{n}, \nabla \mbox{div} \phi> +2 < \nabla \sqrt{n} \otimes \nabla \sqrt{n}, \nabla \phi>.
\end{eqnarray}

Associated to \eqref{eq2}, we can define now the following weak formulation of the momentum equation:
$$  \int_\Omega n_0 \, u_0 \cdot \phi(\cdot,0)\, dx + \int_0^T \int_\Omega (n u\cdot \partial_t\phi + n(u\otimes u): \nabla \phi)  \, dx \, dt
   + \int_0^T \int_\Omega \bigl(p(n)+p_c(n)\bigr){\rm div} \phi \, dx \, dt$$
$$ =  \int_0^T \int_\Omega  \left( 2 \, \eps^2 \sqrt n \,  \nabla \sqrt n \cdot \nabla {\rm div} \phi 
  +  4 \eps^2     \nabla\sqrt n \otimes \nabla\sqrt n : \nabla \phi  +  \int_0^T \int_\Omega  n D(u) : \nabla\phi \right) \, dx \, dt.$$

In this paper, we will first prove that, for a fixed $\varepsilon$, there exists a global weak solution $(n^\varepsilon, u^\varepsilon)$ of \eqref{eq2}. Secondly, we will prove that when $\varepsilon$ tends to zero, $(n^\varepsilon, u^\varepsilon)$ tends to $(n^0,u^0)$ which is a global weak solution of:
\begin{equation} 
\label{eql}
\left\{  
\begin{array} {ll}
\partial_t n^0+\mbox{div}(n^0 \, u^0) = 0, \\[2mm]
\partial_t (n^0 \, u^0)+\mbox{div}(n^0 u^0 \otimes u^0) +\nabla (p(n^0)+  p_c(n^0)) = 2 \, \nu \,  \mbox{div}(n^0 \, D(u^0)), \\[2mm]
n^0|_{t=0} =n_0, \,     \quad  (n^0 \, u^0)|_{t=0} =n_0 \, u_0. \, 
\end{array}
\right.
\end{equation}

Our existence result relies on a careful used of what has been done to prove global existence of weak solutions for the degenerate compressible Navier-Stokes equations by  \cite{BrDe2007} and more recently in the very interesting complete studies \cite{Za2013} (published in \cite{MuPoZa2}). We will then mix these ingredients with some that may be found  in \cite{Jue10}  and \cite{Jiang11} for the construction of approximate solutions. For the asymptotic limit we will strongly use the fact that $\eps$ may vanish letting $\nu$ fixed: the key point being the identity $${\rm div}(n \nabla^{2} \log n)= 2 n \nabla (\Delta\sqrt n/\sqrt n) \, (*) $$ and the presence of a singular pressure close to vacuum allowing to pass to the limit in the nonlinear term $n u\otimes u$. In conclusion that means that the construction of approximate solutions build by A. J\"uengel is consistent with the stability procedure with singular pressure initiated by the works of D. Bresch, B. Desjardins.

\vskip0.2cm
The paper is organized as follow. In Section 2 we state the two main results about existence of solutions and the low Planck constant limit. In Section 3 we prove some a priori estimates and we precise which of them are independent of $\varepsilon$. Section 4 is devoted to the proof of the existence result and will be split in two parts: one for the construction of approximate solution and the second one for the stability of solutions. Finally Section 5 is devoted to the proof of the low Planck constant limit.


\section{Main results}\label{sec.resuts}

In this section we present our two main results. The first one gives the existence of a global weak solution to \eqref{eq2} (in the sense of Definition \ref{def.solfaible1}) without any assumption on $\gamma$ even if the dimension is equal to three. The second one is devoted to the low Planck constant limit and shows that global weak solutions (as defined in Definition \ref{def.solfaible1}) of \eqref{eq2} tends to a global weak solution of \eqref{eql} (in the sense of Definition \ref{def.solfaible2}) when $\varepsilon$ tends to zero. 

Let us first of all give the definitions we will use of weak solution for \eqref{eq2} and \eqref{eql}.
\begin{defi}\label{def.solfaible1}
We say that $(n,u)$ is a weak solution of  \eqref{eq2} if the continuity equation  
\begin{equation}
\left\{  
 \begin{array} {l}
\partial_t n +\mbox{ div } (\sqrt{n}\sqrt{n} u)=0, \\[2mm]
n(0,x)=n_0(x)
\end{array}
\right.
\end{equation}
is satisfied in the sense of distributions and 
the weak formulation of the momentum equation
\begin{multline}\label{eq.solfaible1}
 \int_\Omega n_0 \, u_0  \cdot \phi(\cdot,0)\, dx + \int_0^T \int_\Omega (n u\cdot \partial_t\phi + n(u\otimes u): \nabla \phi) \, dx \, dt
   + \int_0^T \int_\Omega \bigl(p(n)+p_c(n)\bigr){\rm div} \phi \, dx \, dt \\
 =  \int_0^T \int_\Omega \left( 2 \, \eps^2  \,   \sqrt n \,  \nabla \sqrt n \cdot \nabla {\rm div} \phi 
  +  4 \eps^2    \nabla\sqrt n \otimes \nabla\sqrt n : \nabla \phi  + 2 \,   \nu \, n D(u) : \nabla\phi \right) \, dx \, dt.
\end{multline}
holds for any smooth, compactly supported test function $\phi$ such that $\phi(T,.)=0$.
\end{defi}

\begin{defi}\label{def.solfaible2}
We say that $(n^0,u^0)$ is a weak solution of  \eqref{eql} if the continuity equation  
\begin{equation}
\left\{  
 \begin{array} {l}
\partial_t n^0 +\mbox{ div } (\sqrt{n^0}\sqrt{n^0} u^0)=0, \\[2mm]
n^0(0,x)=n_0(x)
\end{array}
\right.
\end{equation}
is satisfied in the sense of distributions and 
the weak formulation of the momentum equation
\begin{multline}\label{eq.solfaible2}
 \int_\Omega n_0 \, u_0 \cdot \phi(\cdot,0)\, dx + \int_0^T \int_\Omega (n^0 u^0\cdot \partial_t\phi + n^0(u^0\otimes u^0): \nabla \phi)   \, dx \, dt \\
   + \int_0^T \int_\Omega \bigl(p(n^0)+p_c(n^0)\bigr){\rm div} \phi   \, dx \, dt
 =  2 \,   \nu \,  \int_0^T \int_\Omega  n^0D(u^0) : \nabla\phi  \, dx \, dt.
\end{multline}
holds for any smooth, compactly supported test function $\phi$ such that $\phi(T,.)=0$.
\end{defi}

Now we can state the main results of this paper. 
First of all, let us introduce the energy of the system which is given by the sum of the kinetic, internal and quantum energies: 
\begin{equation}
\label{def.energie}
E_{\eps}(n,u)= \frac{n}{2}|u|^2+H(n)+H_c(n) + 2  \, \eps ^2 \,  |\nabla \sqrt{n}|^2,
\end{equation}
where $H$ and $H_c$ are given by:
$$H''(n)=\frac{p'(n)}{n} \quad \hbox{and} \quad H_c''(n)=\frac{p_c'(n)}{n}.$$

\begin{remark}
Here, with $p(n)=n^\gamma$ and $\gamma>1$, we have:
$$H(n)=\frac{n^\gamma}{\gamma-1}.$$
\end{remark}
The first main result of this paper  is devoted to the existence of global weak solution in the sense of Definition \ref{def.solfaible1}.

\begin{theorem}\label{thm.existence}
Let $\nu>0, \, \varepsilon >0, \, 1 \leq d \leq 3,  T>0, \gamma \geq 1$. Let $(n_0,u_0)$ such that $n_0 \geq 0$ and $E_\eps(n_0,u_0)<\infty$.
Then there exists a weak solution $(n,u)$ to System \eqref{eq2} in the sense of definition \ref{def.solfaible1} such that
$$n \geq 0 \mbox{ in } \TD^3, \, \sqrt{n} \in L^\infty(0,T;H^1) \cap L^2(0,T;H^2),$$
$$ n \in   L^{\infty}(0,T;L^\gamma),
\,  n^\gamma \in  L^{5/3}(0,T;L^{5/3}),$$
$$\sqrt{n}u \in L^{\infty}(0,T,L^2), \,
 n | \nabla u | \in L^2(0,T;L^2), \,  \sqrt{n} | \nabla u | \in L^2(0,T;L^2),$$
$$\nabla \left( \frac{1}{ \sqrt{n}} \right) \in L^2(0,T;L^2).$$
\end{theorem}

Let us now state the second main result of this paper which gives the convergence of a sequence of global weak solutions $(n^\varepsilon, u^\varepsilon)$ of \eqref{eq2} (in the sense of Definition \ref{def.solfaible1}) to $(n^0,u^0)$ a global weak solution to \eqref{eql} in the sense of Definition \ref{def.solfaible2}.

\begin{theorem} \label{thm.convergence}
Let $1 \leq d \leq 3,~T>0, \, 0< \eps < \nu,  \gamma \geq 1$.
Let $(n_0,u_0)$ such that $n_0 \geq 0$ and $E_\eps(n_0,u_0)<\infty$.
Then for $(n^\eps,u^\eps)$ solution of \eqref{eq2} we have, when $\eps$ tends to $0$:


\begin{align*}
	& (n^{\eps})_\eps \to n^0  \quad \text{strongly in} \ L^2(0,T;L^{\infty}(\Omega)),\\
	& (\sqrt{n^{\varepsilon}})_\eps \rightharpoonup \sqrt{n^0} \quad \text{weakly   in} \ L^2(0,T;H^2(\Omega)),\\
	& (\sqrt{n^{\eps}})_\eps \to \sqrt{n^0} \quad \text{strongly in} \ L^2(0,T;H^{1}(\Omega)), \\
	& (1/\sqrt{n^{\eps}})_\eps \to \frac{1}{\sqrt{n^0}} \quad \text{almost everywhere}, \\ 
	& (\sqrt{n^{\varepsilon}} u ^{\varepsilon})_\eps \to \sqrt{n^0} u^0 \quad \text{strongly in} \ L^2(0,T;L^2(\Omega)),\\
	& (u ^{\varepsilon})_\eps  \rightharpoonup u^0, \quad \text{weakly in} \ L^{p}(0,T;L^{q^\star}(\Omega)),
\end{align*}
with $p={8k}/{(4k+1)},  q^\star= {24k}/{(12k+1)}$ and $(n^0,u^0)$ solution of \eqref{eql}.
\end{theorem}

\begin{remark}
Note that, even if the results are valuable for $1 \leq d \leq 3$, in the proofs we will only focus on the 3d-case. Indeed, the most interesting case in term of difficulties is the 3d-one and the hypothesis $\gamma > 3$ in \cite{Jue10} was only necessary in this case.
\end{remark}


\section{A priori estimates}\label{sec.estimates}

In this section we establish all the a priori estimates we need in order to prove existence and convergence results. We pay a particular attention to the dependence or independence of each one with regards to $\varepsilon$.

Using a formal computation, we easily obtain
\begin{equation}
\label{eq.energie}
\frac{d}{dt} \int_\Omega E_{\eps}(n,u)   \, dx + \nu  \int_\Omega n |D(u)|^2 \, dx =0.
\end{equation}

Directly from \eqref{eq.energie} we deduce the following estimates:
\begin{lemma}\label{lemme1}
Under the hypothesis of Theorem \ref{thm.existence},  there exists a constant $C$ independent of $\varepsilon$, such that    
\begin{eqnarray}
\left\| \sqrt{n}u\right\|_{L^\infty(0,T;L^2(\Omega))}\leq C,\label{eq.lemme1.1}\\
\left\| n^{\gamma} \right\|_{L^\infty(0,T;L^1(\Omega))}\leq C, \label{eq.lemme1.1bis}\\
\left\| \sqrt{n}D(u)\right\|_{L^2(0,T;L^2(\Omega))}\leq C.\label{eq.lemme1.2}
\end{eqnarray}
\end{lemma}

\vskip0.2cm
The following entropy inequality can be shown following the lines of the Bresch-Desjardins entropy and the Bohm potential identity $(*)$.
\begin{proposition}\label{prop.entropy}
Under the hypothesis of Theorem \ref{thm.existence}, we have:
\begin{eqnarray}
\dfrac{d}{dt}
 \int_{\Omega}  \left( \dfrac{n}{2} | u+  \nu \nabla \log n|^2  +H(n)+H_c(n) +  (2 \varepsilon^2+4\nu^{2}) | \nabla  \sqrt{n}|^2  \right) \, dx  \label{eq.entropy}\\
 +  \nu  \int_{\Omega} \left(H''(n) | \nabla n|^2+ H_c''(n) | \nabla n|^2 +  \varepsilon^2 n | \nabla^2 \log n|^2 + 2  n \, |\nabla u |^2\right) \, dx=0.\nonumber
\end{eqnarray}
\end{proposition}

\pr
We have:
$$\partial_t ( n  \nabla \mbox{log} n)+ \mbox{div} (n \, ^t \nabla u)+\mbox{div}(n \, u \otimes \nabla \mbox{log} n )=0.$$
We multiply this new equation by $\nu$, and we add it to:
$$\partial_t (n \, u)+\mbox{div}(n \, u \otimes u) +\nabla (p(n)+p_c(n))-2 \,  \varepsilon^2 \, n \nabla \left( \frac{\Delta \sqrt{n}}{\sqrt{n}} \right)  =  2 \, \nu \,  \mbox{div}(n \, D(u)),$$
to obtain
\begin{eqnarray}
\partial_t(n\, (u+\nu\nabla\log n))&+&\mbox{div}(n\, u \otimes (u+\nu\nabla\log n))\nonumber \\
&+&\nabla(p(n)+p_c(n))-2 \,  \varepsilon^2 \, n \nabla \left( \frac{\Delta \sqrt{n}}{\sqrt{n}} \right) =2 \, \nu \,  \mbox{div}(n \, \nabla u).\nonumber
\end{eqnarray}
Multiplying by $u+\nu \nabla \mbox{log}n$ and integrating over $\Omega$ we have:
\begin{eqnarray}
\int_{\Omega} \partial_{t}\left(n\left(u+\nu \nabla \mbox{log}n\right)\right).\left(u+\nu \nabla \mbox{log}n\right) dx + \int_{\Omega} \mbox{div}\left(nu \otimes \left(u+\nu \nabla \mbox{log}n\right)\right).\left(u+\nu \nabla \mbox{log}n\right) dx \nonumber \\
+\int_{\Omega} \nabla p(n) . \left(u+\nu \nabla \mbox{log}n\right)dx -2\varepsilon^{2} \int_{\Omega} n \nabla \left(\frac{\Delta \sqrt{n}}{\sqrt{n}}\right). \left(u+\nu \nabla \mbox{log}n\right) dx \nonumber \\
=2\nu \int_{\Omega} \mbox{div}\left(n \nabla u\right). \left(u+\nu \nabla \mbox{log}n\right) dx. \label{eq.proof1}
\end{eqnarray}
Moreover, 
\begin{eqnarray}
\hskip-1.6cm\bullet\int_{\Omega} \partial_{t}\left(n\left(u+\nu \nabla \mbox{log}n\right)\right).\left(u+\nu \nabla \mbox{log}n\right) dx&=&\frac{d}{dt} \left(\int_{\Omega} \frac{n}{2} |u+\nu \nabla \mbox{log}n|^{2} dx\right)\nonumber \\
&{}&-\frac{1}{2}\int_{\Omega}\mbox{div}(nu)|u+\nu \nabla \mbox{log}n|^{2} dx \nonumber
\end{eqnarray}
\begin{eqnarray}
\bullet\int_\Omega  \nabla (p(n)+p_c(n)) \cdot (u+2 \nu \nabla \mbox{log}n)dx&= &\int_\Omega   \partial_t (H(n)+H_c(n)) dx \nonumber\\
&{}&+  \nu  \int_\Omega H''(n) | \nabla n|^2 dx +\nu  \int_\Omega H_c''(n) | \nabla n|^2dx \nonumber
\end{eqnarray}
\begin{eqnarray*}
\hskip-0.4cm\bullet- 2 \,  \varepsilon^2 \int_\Omega  n \nabla \left( \frac{\Delta \sqrt{n}}{\sqrt{n}}\right) . (u+2 \nu \nabla \mbox{log}n)
= 2 \varepsilon^2 \int_\Omega   \partial_t((\nabla \sqrt{n})^2) + \nu  \varepsilon^2 \int_\Omega n | \nabla^2 \mbox{log}n|^2 
\end{eqnarray*}
\begin{eqnarray*}
\hskip-0.7cm\bullet\int_{\Omega} \mbox{div}\left(nu \otimes \left(u+\nu \nabla \mbox{log}n\right)\right).\left(u+\nu \nabla \mbox{log}n\right) dx=\frac{1}{2}\int_{\Omega}\mbox{div}(nu)|u+\nu \nabla \mbox{log}n|^{2} dx.
\end{eqnarray*}
Finally, using integration by parts,
\begin{eqnarray*}
\bullet \ 2\nu \int_{\Omega} \mbox{div}\left(n\nabla u\right). \left(u+\nu \nabla \mbox{log}n\right) dx & = &
-2\nu\int_{\Omega} n \left|\nabla u\right|^{2}dx + 2\nu^{2} \int_{\Omega} u . \mbox{div}\left(n \nabla^{2} \mbox{log} n\right)dx\\
 &= & -2\nu\int_{\Omega} n \left|\nabla u\right|^{2}dx + 4\nu^{2} \int_{\Omega} nu \nabla\left(\frac{\Delta \sqrt{n}}{\sqrt{n}}\right)dx \\
 &= & -2\nu\int_{\Omega} n \left|\nabla u\right|^{2}dx - 4\nu^{2} \int_{\Omega} \mbox{div} (nu) \left(\frac{\Delta \sqrt{n}}{\sqrt{n}}\right)dx \\
 & = & -2\nu\int_{\Omega} n \left|\nabla u\right|^{2}dx + 4\nu^{2} \int_{\Omega} \partial_{t}n \left(\frac{\Delta \sqrt{n}}{\sqrt{n}}\right)dx \\
& = & -2\nu\int_{\Omega} n \left|\nabla u\right|^{2}dx - 4\nu^{2} \int_{\Omega} \partial_{t} \left(|\nabla \sqrt{n}|^{2}\right)dx.
\end{eqnarray*}
Using all the above inequalities in \eqref{eq.proof1} we obtain proposition \ref{prop.entropy}.
 \cqfd

Directly from Proposition \ref{prop.entropy} we deduce the following estimates:
\begin{lemma}\label{lemme2}
Under the hypothesis of Theorem \ref{thm.existence},  there exists a constant $C$ independent of $\varepsilon$, such that    
\begin{eqnarray}
\left\| \sqrt{n}(u+\nu \nabla \log n)\right\|_{L^\infty(0,T;L^2(\Omega))}\leq C,\label{eq.lemme2.1}\\
\left\| \nabla \sqrt{n}\right\|_{L^\infty(0,T;L^2(\Omega))}\leq C,\label{eq.lemme2.2}\\
\left\| \nabla n^{\gamma/2}  \right\|_{L^2(0,T;L^2(\Omega))}\leq C, \label{autre} \\
\left\| \sqrt{n} \nabla u\right\|_{L^2(0,T;L^2(\Omega))}\leq C, \label{eq.lemme2.3} \\
\varepsilon \left\| \sqrt{n} \nabla^2 \log n\right\|_{L^2(0,T;L^2(\Omega))}\leq C, \label{eq.lemme2.4} \\
\varepsilon \left\| \nabla^2 \sqrt{n}\right\|_{L^2(0,T;L^2(\Omega))}\leq C. \label{eq.lemme2.5}
\end{eqnarray}
\end{lemma}

\begin{remark}
Contrary to \eqref{eq.energie}, the entropy equality \eqref{eq.entropy} allows us to obtain \eqref{eq.lemme2.2} with $C$ independent of $\varepsilon$ thanks to the coefficient $(2 \varepsilon^2+4\nu^{2})$ which can be minor by $4\nu^{2}$.
\end{remark}

\pr
Estimates \eqref{eq.lemme2.1}, \eqref{eq.lemme2.2}, \eqref{eq.lemme2.3} and \eqref{eq.lemme2.4} are just direct consequence from \eqref{eq.entropy}.
Estimate \eqref{autre} derived from 
$$H''(n) |\nabla n|^2=\frac{p'(n)|\nabla n|^2}{n}  = \gamma n^{\gamma-2}|\nabla n|^2 \quad \text{and} \quad \nabla n^{\gamma/2}=\frac{\gamma}{2} n^{\gamma/2-1} \nabla n.$$

 The last one is obtain using the following identity (which may be found for instance in \cite{Jue10})
$$\int_0^T \int_{\tore} n |\nabla^2 \mbox{log} n|^2 \,  dx  \, dt \geq \int_0^T \int_{\tore} | \nabla^2 \sqrt{n}|^2 \, dx  \, dt$$
and \eqref{eq.lemme2.4}.
\cqfd

\vskip0.2cm
Let us now show some estimates concerning the two different pressures.
\begin{lemma}\label{lemme3}
Under the hypothesis of Theorem \ref{thm.existence}, assuming that  $d=3$, there exists a constant $C$ independent of $\varepsilon$, such that    
\begin{eqnarray}
\left\| n^\gamma \right\|_{L^{5/3}(0,T;L^{5/3}(\Omega))}\leq C,\label{eq.lemme3.1}\\
\left\|p_c(n)\right\|_{L^{5/3}(0,T;L^{5/3}(\Omega))}\leq C.\label{eq.lemme3.2}
\end{eqnarray}
\end{lemma}

\pr
Using \eqref{eq.lemme1.1bis} and \eqref{autre} 
and Sobolev embedding for $d=3$, we obtain
$$\left\| n^{\gamma/2} \right\|_{L^{2}(0,T;L^{6}(\Omega))}\leq C,$$
which can also be written
$$\left\| n^\gamma \right\|_{L^{1}(0,T;L^{3}(\Omega))}\leq C.$$
Using interpolation, last inequality and \eqref{eq.lemme1.1bis} give
$$\left\| n^\gamma \right\|_{L^{5/3}(0,T;L^{5/3}(\Omega))}\leq C.$$
Let us recall that   
$$p_c(n)=\frac{c}{-4k}n^{-4k}, \, n \leq 1, \,  \, k>1.$$
Let $\zeta$ be  a smooth function such that  
$$\zeta(y)=y \mbox{ for } y \leq 1/2 \mbox{  and } \zeta(y)=0 \mbox{ for }y>1.$$
Following the idea developed in \cite{BrDe2007} (page 71) and more recently in \cite{Za2013} (page 62), using \eqref{eq.energie} and \eqref{eq.entropy} we obtain (through $H_c$) for a constant $C$ independent of $\varepsilon$
$$\int_0^T \int_{\tore} | \nabla \zeta(n)^{-2k}|^2 \, dx \, dt\leq C.$$
That means that $\nabla \zeta(n)^{-2k} \in L^2(0,T;L^2(\Omega))$ and close to vacuum, 
\begin{equation}
\label{eq.proof.0}
\sup_t \int_{\tore} n^{-4k} \, dx  \leq C,
\end{equation}
with $C$ independent of $\varepsilon$. This gives $\zeta (n)^{-2k} \in L^{\infty}(0,T;L^2(\Omega))$, and therefore we have $\zeta(n)^{-2k} \in L^2(0,T;H^1(\Omega))$.

Using Sobolev embedding in dimension $d=3$, 
\begin{equation}
\label{eq.proof.0bis} 
\|\zeta(n)^{-2k}\|_{L^2(0,T;L^6(\Omega))}\leq C,
\end{equation}
with still $C$ independent of $\varepsilon$. Now using as previously interpolation and \eqref{eq.proof.0}, \eqref{eq.proof.0bis} we obtain \eqref{eq.lemme3.2}. 
This ends the proof of Lemma \ref{lemme3}.
\cqfd

Let us now prove some estimates that will be used to pass to the limit in the nonlinear term $n u\otimes u$. That means to get the strong convergence of $\sqrt n u$ in  $L^2(0,T;L^2(\Omega))$. Such convergence is required for the stability process to build global weak solution but also for the asymptotic analysis $\eps$ tends to $0$.  This follows the lines that may be found in \cite{BrDe2007} and \cite{Za2013} but taking care of the new quantum term. We will only recall the main  steps and refer the interested reader to very nice PhD Thesis  \cite{Za2013} for all details.

The idea is to prove that $\sqrt n u $ is uniformly bounded in $L^{p'}(0,T,L^{q'}(\Omega))$ for $p',q'> 2$. This will be used with almost pointwise convergence to deduce the strong convergence in $L^2(0,T;L^2(\Omega))$. The main steps are the following: look at the estimates on $\nabla u$ (using the estimates on $\sqrt n \nabla u$ and the estimates on some negative power of $n$) and therefore on $u$ by Sobolev embeddings. Mix it with the bounds on $\sqrt n u$ and $n$  to increase the uniform integrability in space on $\sqrt n u$ uniformly with respect to $\eps$.

\begin{lemma}\label{lemme4}
Under the hypothesis of Theorem \ref{thm.existence}, assuming that $d=3$, there exists a constant $C$ independent of $\varepsilon$, such that    
\begin{equation}
\label{eq.lemme4}
\left\| \nabla u \right\|_{L^{p}(0,T;L^{q}(\Omega))}\leq C,
\end{equation}
with $p=8k/(4k+1)$ and $q=24k/(12k+1)$.
\end{lemma}

\pr
Let us write 
$$\nabla u=\frac{1}{\sqrt{n}}\sqrt{n}\nabla u.$$
Using \eqref{eq.lemme2.3} we have an estimate in $L^2(0,T;L^2(\Omega))$ for $\sqrt{n}\nabla u$. Then it remains to obtain an estimate for $1/\sqrt{n}$ in the appropriate space. 
%
Using \eqref{eq.proof.0bis} and the remark that close to vacuum,
$$ \|\zeta(n)^{-2k}\|_{L^2(0,T;L^6(\Omega))}^2=\int_0^T\left(\int_\Omega \left(\frac{1}{\sqrt{n}}\right)^{24k}dx\right)^{8k/24k}dt,$$
we obtain
\begin{equation}
\label{eq.proof.1}
\| 1/\sqrt{n}\|_{L^{8k}(0,T;L^{24k}(\Omega))} \leq C,
\end{equation}
with $C$ independent of $\varepsilon$. 

As previously said, since $\nabla u=(1/\sqrt{n})\sqrt{n}\nabla u$, using \eqref{eq.lemme2.3} and \eqref{eq.proof.1}, we have the result.
\cqfd

Using Sobolev embedding, a direct consequence of Lemma \ref{lemme4} is the following one.
\begin{lemma}\label{lemme5}
Under the hypothesis of Theorem \ref{thm.existence}, assuming that $d=3$, there exists a constant $C$ independent of $\varepsilon$, such that    
\begin{equation}
\label{eq.lemme5}
\left\| u \right\|_{L^{p}(0,T;L^{q^\star}(\Omega))}\leq C,
\end{equation}
with $p=8k/(4k+1)$ and $q^\star=24k/(12k+1)$.
\end{lemma}

\vskip0.2cm
In the following lemma, we have  an estimate on the gradient of negative power of $n$:

\begin{lemma}\label{lemme6}
Under the hypothesis of Theorem \ref{thm.existence}, there exists a constant $C$ independent of $\varepsilon$, such that    
\begin{equation}
\label{eq.lemme6}
\left\| \nabla (1/\sqrt{n}) \right\|_{L^{2}(0,T;L^{2}(\Omega))}\leq C.
\end{equation}
\end{lemma}

\vskip0.2cm
\pr
We know from  \eqref{eq.lemme2.2} that $\nabla \sqrt{n}  \in L^2(\Omega)$.
But $$\nabla \left( \frac{1}{\sqrt{n}} \right)=-\frac{1}{2} \frac{\nabla n} {n^{3/2}}= - \frac{\nabla \sqrt{n}}{n}.$$
Then, if $n >1$, $\frac{1}{n} <1$ and $ \nabla (1/\sqrt{n}) \in L^2(\Omega)$.

Now look at the case $n \leq 1$.
From Proposition \ref{prop.entropy}, we have 
$\int_0^T \int_{\Omega} H''_c(n) |\nabla n |^2 \, dx < + \infty$ where $H''_c(n) = \frac{p'_c(n)}{n}$.
For $n \leq 1$, $p'_c(n)= c n ^{-4k-1}$
then  $H''_c(n) = c n^{-4k-2}$, $$H''(c)  |\nabla n |^2= c \left| \frac{\nabla n}{n^{2k+1}}\right|^2= \frac{c}{4k^2}\left |\nabla \left(\frac{1}{n^{2k}}\right)\right|^2.$$
Make the connection between $\nabla \left(\frac{1}{n^{2k}}\right)$ and $\nabla \left(\frac{1}{\sqrt{n}}\right)$ :
\begin{eqnarray*}
\nabla \left(\frac{1}{\sqrt{n}}\right) & = & \nabla\left(\frac{1}{n^{2k}} n^{2k-1/2}\right)  \\
 &= & n^{2k-1/2}  \nabla\left(\frac{1}{n^{2k}}\right) + \frac{1}{n^{2k}} \nabla\left(n^{2k-1/2}\right) \\
& = & n^{2k-1/2}  \nabla\left(\frac{1}{n^{2k}}\right) +(2k-1/2) n^{-2k} \nabla(n) n^{2k-3/2} \\
 & = & n^{2k-1/2}  \nabla\left(\frac{1}{n^{2k}}\right) +(2k-1/2)  \nabla(n) n^{-3/2}=n^{2k-1/2}  \nabla\left(\frac{1}{n^{2k}}\right) -2(2k-1/2)  \nabla\left(\frac{1}{\sqrt{n}}\right),
\end{eqnarray*}
we have
$$(1+4k-1)\nabla \left(\frac{1}{\sqrt{n}}\right)= n^{2k-1/2}  \nabla\left(\frac{1}{n^{2k}}\right),$$
and 
$$\left|\nabla \left(\frac{1}{\sqrt{n}}\right)\right|^2= \frac{1}{16k^2} n^{4k-1} \left| \nabla\left(\frac{1}{n^{2k}}\right)\right|^2=\frac{1}{16k^2} n^{4k-1} 4k^2 \frac{1}{c} H''_c(n) |\nabla n |^2=\frac{1}{4c} n^{4k-1}H''_c(n) |\nabla n |^2.$$
Such as $n \leq 1$ and $\int_0^T\int_{\Omega} H''_c(n) |\nabla n |^2 \, dx < + \infty$, we have $\nabla \left(\frac{1}{\sqrt{n}}\right) \in L^2(0,T;L^2(\Omega))$.
 \cqfd

\vskip0.2cm
Using previous lemma, we are now able to establish the following proposition.

\begin{proposition}\label{prop2}
Under the hypothesis of Theorem \ref{thm.existence}, assuming that $d=3$, there exists a constant $C$ independent of $\varepsilon$, such that    
\begin{equation}
\label{eq.prop.2}
\|\sqrt{n}u\|_{L^{p'}(0,T;L^{q'}(\Omega))}\leq C,
\end{equation}
with $p', q' >2$.
\end{proposition}

\pr
Let $r>0$ to be chosen later on. We write 
$$\sqrt{n}u=\left(\sqrt{n}u\right)^{2r}u^{1-2r}n^{1/2-r}.$$
Using \eqref{eq.lemme1.1bis}
$$\| n^{1/2-r}\|_{L^\infty(0,T;L^{\gamma/(1/2-r)}(\Omega))} \leq C, $$
with $C$ a constant independent of $\varepsilon$. Using \eqref{eq.lemme1.1}, \eqref{eq.lemme5} and \eqref{eq.lemme6}, choosing
$$ \frac{1}{p'}=\frac{1-2r}{p}\, \hbox{and} \, \frac{1}{q'}=\frac{2r}{2}+\frac{1-2r}{q^\star}+\frac{1/2-r}{\gamma},$$
with $r=2/3$, we get the conclusion. Indeed, with such definition of $p'$ and $q'$, the condition $q'>2$ is equivalent to $r> 1/2$ and the condition $p'>2$ is equivalent to $r>1/(8k+2)$ with $1/(8k+2)<1/10$ since $k>1$.
\cqfd

\begin{remark}
Estimate \eqref{eq.prop.2} is the one proved in \cite{Za2013}. Note that it is also possible to choose
drag terms of form $r_0 u +r_1 n |u|u$ instead of cold pressure. Such a choice provides a bound on 
 $n \, u^{2+\delta}$ in $L^1 (0,T; L^1(\Omega))$ uniformly with $\delta >0$ directly from the energy estimate (cf the trick from  Mellet and Vasseur  \cite{MV}).
\end{remark}


\section{Proof of global existence of solutions}\label{sec.existence}

In this section we prove existence of a global weak solution to \eqref{eq2} {\it i.e.} we prove Theorem \ref{thm.existence}. The proof can be split into two parts: in the first one we show how to construct an approximate solution and in the second one we show the stability of solution. This is the classical way to prove existence of solution for Navier-Stokes equations. For clarity of presentation we will detail each part in one subsection.


\subsection{Construction of an approximate solution}\label{subsec.construction}

In this subsection we present the construction of an approximate solution which follows the lines given in \cite{Jue10}. In few words, the approximate solution is build through a fixed point argument at the level of the Galerkin approximate system, and, thanks to appropriate uniform estimates with respect to the Galerkin parameters, we can pass to the limit and obtain existence of global weak solution. 

\vskip0.2cm
Let us now describe precisely the procedure. As in \cite{Jue10}, let us first transform \eqref{eq2} by the use of the so-called effective velocity
$$w=u+\nu \nabla \log n. $$
Then a simple computation shows that \eqref{eq2} is equivalent to: for $x \in \Omega$ and $t>0$:
\begin{equation} 
\label{eqw}
\left\{  
 \begin{array}{ll}
\partial_t n+\mbox{div}(n \, w) = \nu \Delta n, \\[2mm]
\partial_t(n \, w)+\mbox{div}(n \, w \otimes w) +\nabla_x (p(n) + p_c(n))
        -2 \,    \varepsilon_{0} \, n \nabla \left( \dis  \frac{\Delta \sqrt{n}}{\sqrt{n}} \right)  =  \nu \,  \Delta(nw),Ê\\[2mm]
n|_{t=0} =n_0, \,  \quad (n \, w)|_{t=0} =n_{0}w_{0}, \, 
\end{array}
\right.
\end{equation}
with $w_{0}=u_{0}+\nu \nabla \log n_{0}$ and $ \varepsilon_{0}=\varepsilon^{2}-\nu^{2}$. We associate to \eqref{eqw} the following definition of weak solution.
\begin{defi}\label{def.solfaible3}
We say that $(n,w)$ is a weak solution of  \eqref{eqw} if the equation  
\begin{equation}
\label{eq.def.solfaible3}
\left\{  
 \begin{array} {l}
\partial_t n +\mbox{ div } (n w)=\nu \Delta n, \\
n(0,x)=n_0(x)
\end{array}
\right.
\end{equation}
is satisfied in the sense of distributions and 
the weak formulation 
\begin{multline}\label{eq.solfaible3}
 \int_\Omega n_0w_{0}\cdot \phi(\cdot,0)\, dx + \int_0^T \int_\Omega (n w\cdot \partial_t\phi + n(w\otimes w): \nabla \phi)  \, dx \, dt  \\
   + \int_0^T \int_\Omega \bigl(p(n)+p_c(n)\bigr){\rm div} \phi  \, dx \, dt 
 =   \,   \int_0^T \int_\Omega  \left( 2 \, \eps_{0}^2 \frac{\Delta\sqrt n}{\sqrt{n}} \,  {\rm div} (n\,\phi) 
  +  \nu \,   \nabla(nw) : \nabla\phi \right) \, dx \, dt.
\end{multline}
holds for any smooth, compactly supported test function $\phi$ such that $\phi(T,.)=0$.
\end{defi}
Now the goal is to prove existence of solution to \eqref{eqw} in the sense of Definition \ref{def.solfaible3}. To this end we use the same technic as the one described Sections 3, 4, 5 and 6 in \cite{Jue10}.

\subsubsection{Global existence of solutions}

\noindent

Let $T>0$ and $(e_{p})$ an orthonormal basis of $L^{2}(\Omega)$. Note that $(e_{p})$ is also an orthogonal basis of $H^{1}(\Omega)$. We introduce the finite space $X_{N}=span\{e_{1}, \cdots, e_{N}\}$, for $N \in \mathbb{N}^*$. Let $(n_{0},w_{0}) \in \mathcal{C}^{\infty}(\Omega)^{2}$ some initial data such that $n_{0}(x)\geq \delta >0$ for all $x \in \Omega$ and for some $\delta >0$. Let $v \in \mathcal{C}^{0}([0,T];X_{N})$ a given velocity, $v$ can be written for $x \in \Omega$ and $t\in [0,T]$
$$v(x,t)=\sum_{i=1}^{N} \lambda_{i}(t)e_{i}(x),$$
for some functions $\lambda_{i}$. The norm of $v$ in $\mathcal{C}^{0}([0,T];X_{N})$ can be formulated as
$$\|v \|_{\mathcal{C}^{0}([0,T];X_{N})} =\max_{t \in [0,T]} \sum_{i=1}^{N} |\lambda_{i}(t)|,$$
which has for consequence that $v$ is bounded in $\mathcal{C}^{0}([0,T];\mathcal{C}^{n}(\Omega))$ for any $n \in \mathbb{N}$, and there exists a constant $C$ depending on $n$ such that
\begin{equation}
\label{eq.exist.1}
\|v \|_{\mathcal{C}^{0}([0,T];\mathcal{C}^{n}(\Omega))}\leq C \|v \|_{\mathcal{C}^{0}([0,T];L^{2}(\Omega))}.
\end{equation}
As in \cite{Jue10}, the approximate system is defined as follows. Let $n \in \mathcal{C}^{1}([0,T];\mathcal{C}^{3}(\Omega))$ be the classical solution to
\begin{equation}
\label{eq.exist.2}
\left\{  
 \begin{array} {l}
\partial_t n +\mbox{ div } (n v)=\nu \Delta n, x \in \Omega, t \in [0,T]\\
n(0,x)=n_0(x), x \in \Omega.
\end{array}
\right.
\end{equation}
Using the maximum principle, which provides lower and upper bounds, the assumption $$n_{0}\geq \delta >0$$ and \eqref{eq.exist.1}, $n$ is strictly positive and for all $(x,t) \in \Omega \times [0,T]$,
$$0<\underline{n}(c)\leq n(x,t) \leq \overline{n}(c). $$
Moreover, for a given $n_N$ solution of \eqref{eq.exist.2}, we are looking for a function $w_{N} \in \mathcal{C}^{0}([0,T];X_{N})$ such that
\begin{multline}\label{eq.exist.3}
- \int_\Omega n_0w_{0}\cdot \phi(\cdot,0)\, dx = \int_0^T \int_\Omega \biggl(n_N w_{N}\cdot \partial_t\phi + n_N(v\otimes w_{N}): \nabla \phi + (p(n_N)+p_c(n_N)){\rm div} \phi  \\
 - 2 \, \eps_{0}  \, \frac{\Delta\sqrt n_N}{\sqrt{n_N}} \,  {\rm div} (n_N\,\phi) -  \nu \, \nabla(n_Nw_{N}) : \nabla\phi-\delta(\nabla w_{N}: \nabla \phi + w_{N} \cdot \phi)\biggl) dx \, dt.
\end{multline}
As detailed in \cite{Jue10}, using a Banach fixed point theorem, there exists a unique local-in-time solution $(n_{N},w_{N})$ to \eqref{eq.exist.2} and \eqref{eq.exist.3} with $n_{N} \in \mathcal{C}^{1}([0,T'];\mathcal{C}^{3}(\Omega))$ and $w_{N} \in \mathcal{C}^{1}([0,T'];X_{N})$, for $T' \leq T$. In order to prove the global nature of the solution $(n_{N},w_{N})$ constructed above we use the following energy estimate.
\begin{lemma}\label{lemme.exist.1}
Let $T' \leq T$, and let $n_{N} \in \mathcal{C}^{1}([0,T'];\mathcal{C}^{3}(\Omega)),\, w_{N} \in \mathcal{C}^{1}([0,T'];X_{N})$ be a local-in-time solution to \eqref{eq.exist.2} and \eqref{eq.exist.3} with $n=n_{N}$ and $v=w_{N}$. Then
\begin{multline}
\label{eq.lemmeexist}
\frac{dE_{\eps_{0}}}{dt}(n_{N},w_{N})+ \nu \int_{\Omega} \left(n_{N}|\nabla w_{N}|^{2} + H''(n_{N})|\nabla n_{N}|^{2}+H_{c}''(n_{N})|\nabla n_{N}|^{2}\right) dx \\
+\eps_{0} \, \nu \int_{\Omega} n_{N}|\nabla^{2}\log n_{N}|^{2} dx 
+\delta \int_{\Omega} (|\nabla w_{N}|^{2}+|w_{N}|^{2})dx=0,
\end{multline}
where
$$E_{\eps_{0}}(n_N,w_N)=\int_\Omega \left(\frac{n_N}{2}|w_N|^2+H(n_N)+H_c(n_N)+2\eps_0^2|\nabla \sqrt{n_N}|^2\right)dx.$$
\end{lemma}
We refer the interested reader to \cite{Jue10} for the proof of such lemma.
Then the both limits $N \rightarrow \infty$ and $\delta \rightarrow 0$ are considered separately. 

\begin{remark}
Note that as in \cite{Jue10}, equations for $u$ and $w$ being equivalent, using \eqref{eq.lemmeexist}, the estimates obtained previously on $u$ are still true for $w_N$ and $w_\delta$ with constants $C$ independent of $N$ and $\delta$.
\end{remark}

\subsubsection{Limits $N \rightarrow \infty$ and $\delta \rightarrow 0$}

First we perform the limit $N \rightarrow \infty$, $\delta>0$ being fixed. This is achieved by the use of regularities of the solution and Aubin-Simon's lemma. Here, the use of a cold pressure term avoid the hypothesis $\gamma>3$ for $d=3$ which was required in \cite{Jue10}.

\vskip0.2cm
We can prove the following proposition.
\begin{proposition}\label{prop.cvgce.N}
Under the hypothesis of Theorem \ref{thm.existence}, for a fixed $\eps$, up to a subsequence, the following convergences hold when $N$ tends to $\infty$.
\begin{align*}
&\sqrt{n_{N}} \to \sqrt{n_\delta}, \quad  \hbox{strongly in } \ L^{2}(0,T;H^{1}(\Omega)),\\
&p(n_{N}) \to p(n_\delta), \quad \hbox{strongly in } \ L^{1}(0,T;L^{1}(\Omega)), \\
&p_{c}(n_{N}) \to p_{c}(n_\delta), \quad \hbox{strongly in } \ L^{1}(0,T;L^{1}(\Omega)),\\
&1/\sqrt{n_{N}} \to 1/\sqrt{n_\delta}, \quad \hbox{almost everywhere }, \\
&\sqrt{n_{N}}\,w_{N} \to \sqrt{n_\delta}\,w_\delta, \quad  \hbox{strongly in } \ L^{2}(0,T;L^{2}(\Omega)), \\
&\nabla w_N \rightharpoonup \nabla w_\delta, \quad \hbox{weakly in } \ L^{p}(0,T;L^{q}(\Omega)),\\
&w_N \rightharpoonup w_\delta, \quad \hbox{weakly in } \ L^{p}(0,T;L^{q^\star}(\Omega)),
\end{align*}
with $p,\, q$ and $q^\star$ given in Lemmas \ref{lemme4} and \ref{lemme5}.
\end{proposition}

\pr
First of all, let say that using \eqref{eq.lemmeexist} and the same technics as used in Section 3, we can prove estimates \eqref{eq.lemme1.1}, \eqref{eq.lemme2.2}-\eqref{eq.lemme5} with $n=n_N$ and $u=w_N$ and with constants $C$ which are all independent of $N$ and $\delta$ but can depend on $\eps$.

Using \eqref{eq.lemme1.1}, \eqref{eq.lemme2.2}, \eqref{eq.lemme2.3}, \eqref{eq.lemme2.5}, \eqref{eq.proof.1} and rewriting equation \eqref{eq.def.solfaible3} as
$$\partial_t (\sqrt{n_N})+\frac{1}{2\sqrt{n_N}}\hbox{div}(n_N\,w_N)=\nu \left(\Delta \sqrt{n_N} + \frac{|\nabla \sqrt{n_N}|^2}{\sqrt{n_N}}\right) ,$$ 
we can show that
\begin{equation}
\label{proof.exist.2}
\| \partial_{t}(\sqrt{n_{N}})\|_{L^{2}(0,T;H^{-1}(\Omega))}\leq C,
\end{equation}
with $C$ a constant independent of $\eps$. Then using \eqref{eq.lemme2.5} and \eqref{proof.exist.2} we can apply the Aubin-Simon's Lemma (see \cite{Simon}) to obtain the strong convergence of $\sqrt{n_{N}}$ to $\sqrt{n}$ in $L^{2}(0,T;H^{1}(\Omega))$. Note that \eqref{eq.lemme2.5} being badly dependent on $\eps$, this strong convergence is only true for a fixed $\eps$.

Using \eqref{eq.lemme3.1}, \eqref{eq.lemme3.2} and the almost everywhere convergence of $n_{N}$ (which is a direct consequence of the previous strong convergence for $\sqrt{n_{N}}$), we obtain the strong convergence of the two pressures in $L^{1}(0,T;L^{1}(\Omega))$.

Rewriting equation \eqref{eq.def.solfaible3} as
$$ \partial_t \left( \frac{1}{\sqrt{n_N}} \right)+ \nabla\left(\frac{w_N}{\sqrt{n_N}}\right)+\frac{3}{2\sqrt{n_N}}\hbox{div}(w_N)=-\nu \left(\frac{\Delta \sqrt{n_N}}{n_N} + \frac{|\nabla \sqrt{n_N}|^2}{{n_N}^{3/2}}\right) ,$$
and using \eqref{eq.lemme1.1},  \eqref{eq.lemme2.2}, \eqref{eq.lemme2.3}, \eqref{eq.lemme2.5}, \eqref{eq.proof.0} and \eqref{eq.proof.1}, we have
\begin{equation}
\label{proof.exist.3}
|| \partial_t  \left( \frac{1}{\sqrt{n_{N}}} \right) ||_{L^{\infty}(0,T;W^{-1,1}(\Omega))}\leq C,
\end{equation}
with $C$ a constant independent of $N$ and $\delta$ but which can depend on $\eps$. Then \eqref{proof.exist.3} and \eqref{eq.proof.1} allow us to apply the Aubin-Simon's Lemma and to obtain the almost everywhere convergence of $1/\sqrt{n_N}$.

Moreover, we have 
$$\nabla(n_{N}w_{N})=n_{N}\nabla w_{N}+ w_{N} \nabla n_{N}=\sqrt{n_{N}} \sqrt{n_{N}}\nabla w_{N}+2 \, \sqrt{n_{N}} \,  w_{N} \nabla (\sqrt{n_{N}}) \in L^2(0,T; L^1(\Omega)).$$
Therefore $n_{N}w_{N} \in L^2(0,T;W^{1,1}(\Omega))$. Using the second equation in \eqref{eqw}, we can get
information on $\partial_t(n_{N}w_{N})$ and therefore through Aubin-Lions-Simon's Lemma convergence almost everywhere of $n_{N}w_{N}$. With \eqref{eq.prop.2} and since $\sqrt{n_{N}} w_{N}$ converges almost everywhere then $(\sqrt{n_{N}} w_{N})^2$ converges strongly in $L^1(0,T;L^1(\Omega))$ and therefore  $\sqrt{n_{N}} w_{N}$ converges strongly in  $L^2(0,T;L^2(\Omega))$.

Finally the last two weak convergences directly come from estimates \eqref{eq.lemme4} and \eqref{eq.lemme5}.
\cqfd

Using Proposition \ref{prop.cvgce.N}, passing to the limit $N$ tends to $\infty$, we obtain the existence of a solution $(n_\delta,w_\delta)$ satisfying \eqref{eq.lemmeexist}. Note that since estimates \eqref{eq.lemme1.1}, \eqref{eq.lemme2.2}-\eqref{eq.lemme5} for $n=n_N$ and $u=w_N$ were true for constants $C$ independent of $N$ and $\delta$, they are still true for $n=n_\delta$ and $u=w_\delta$ and with constants $C$ independent of $\delta$ but which can depend on $\eps$. In a same way that previously this allows us to obtain the following convergences.

\begin{proposition}\label{prop.cvgce.delta}
Under the hypothesis of Theorem \ref{thm.existence}, for a fixed $\eps$, up to a subsequence, the following convergences hold when  $\delta$ tends to $0$.
\begin{align*}
&\sqrt{n_{\delta}} \to \sqrt{n}, \quad \text{strongly in } \ L^{2}(0,T;H^{1}(\Omega)),\\
&p(n_{\delta}) \to p(n), \quad \text{strongly in } \ L^{1}(0,T;L^{1}(\Omega)), \\
&p_{c}(n_{\delta}) \to p_{c}(n), \quad \text{strongly in } \ L^{1}(0,T;L^{1}(\Omega)), \\
&1/\sqrt{n_{\delta}} \to 1/\sqrt{n}, \quad \text{almost everywhere },\\
&\sqrt{n_{\delta}}\,w_{\delta} \to \sqrt{n}\,w, \quad \text{strongly in } \ L^{2}(0,T;L^{2}(\Omega)). \\
&\nabla w_\delta \rightharpoonup \nabla w, \quad \text{weakly in } \ L^{p}(0,T;L^{q}(\Omega)),\\
&w_\delta \rightharpoonup w, \quad \text{weakly in } \ L^{p}(0,T;L^{q^\star}(\Omega)),
\end{align*}
with $p,\, q$ and $q^\star$ given in Lemmas \ref{lemme4} and \ref{lemme5}.
\end{proposition}

The proof of this proposition being really similar to the one for Proposition \ref{prop.cvgce.N}, we skip it here. Thanks to Proposition \ref{prop.cvgce.delta} we obtain the existence of solution to \eqref{eqw}  in the sense of Definition \ref{def.solfaible3} without the hypothesis $\gamma>3$ even if we have used the same procedure as the one presented in \cite{Jue10}. This achieves the construction of an approximate solution using the relation between $u$ and $w$.


\subsection{Stability}\label{subsec.stability}

In this subsection we look at the stability. Then let us assume that there exists a sequence $(n_\tau,u_\tau)$ of global weak solutions to \eqref{eq2} satisfying uniformly the energy and entropy inequalities. The goal is to prove that there exists a subsequence which converges to a global weak solution of \eqref{eq2} satisfying also the energy and entropy inequalities. Then it remains to prove that we can pass to the limit in 
$$n_\tau u_\tau, \, n_\tau u_\tau \otimes u_\tau,  \, n_\tau D(u_\tau), \, p(n_\tau), \, p_c(n_\tau), \, \nabla\sqrt n_\tau \otimes \nabla\sqrt n_\tau$$
and in the quantum term $n_\tau \nabla (\displaystyle \frac{1}{\sqrt{n_\tau}} \Delta \sqrt{n_{\tau}})$ namely  in the terms with $ \sqrt{n_\tau}$ or $\nabla \sqrt{n_\tau}$.
To this end we need some strong convergences which can be obtained by the use of Aubin-Simon's  Lemma. More precisely, we have the following proposition.
\begin{proposition}\label{prop3}
Under the hypothesis of Theoerem \ref{thm.existence}, we have for a fixed $\varepsilon$
\begin{align*}
&\sqrt{n_{\tau}} \to \sqrt{n}, \quad \hbox{strongly in} \ L^{2}(0,T;H^{1}(\Omega)),\\
&p(n_{\tau}) \to p(n), \quad \hbox{strongly in} \ L^{1}(0,T;L^{1}(\Omega)),\\
&p_{c}(n_{\tau}) \to p_{c}(n), \quad \hbox{strongly in} \ L^{1}(0,T;L^{1}(\Omega)), \\
&1/\sqrt{n_{\tau}} \to 1/\sqrt{n}, \quad \hbox{almost everywhere},\\
&\sqrt{n_{\tau}}\,u_{\tau} \to \sqrt{n}\,u, \quad \hbox{strongly in} \ L^{2}(0,T;L^{2}(\Omega)).
\end{align*}
\end{proposition}

\pr
Using \eqref{eq.lemme1.1}, \eqref{eq.lemme2.3} and the mass equation, we can show that
\begin{equation}
\label{proof.stability.1}
\| \partial_{t}(\sqrt{n_{\tau}})\|_{L^{2}(0,T;H^{-1}(\Omega))}\leq C,
\end{equation}
with $C$ a constant independent of $\eps$. Then using \eqref{eq.lemme2.5} and \eqref{proof.stability.1} we can apply the Aubin-Simon's  Lemma (see \cite{Simon}) to obtain the strong convergence of $\sqrt{n_{\tau}}$ to $\sqrt{n}$ in $L^{2}(0,T;H^{1}(\Omega))$. Note that \eqref{eq.lemme2.5} being badly dependent on $\eps$, this strong convergence is only true for a fixed $\eps$.

Using \eqref{eq.lemme3.1}, \eqref{eq.lemme3.2} and the almost everywhere convergence of $n_{\tau}$ (which is a direct consequence of the previous strong convergence for $\sqrt{n_{\tau}}$), we obtain the strong convergence of the two pressures in $L^{1}(0,T;L^{1}(\Omega))$.

Again using the mass equation and \eqref{eq.lemme1.1}, \eqref{eq.lemme2.3} and \eqref{eq.proof.0}, we have
\begin{equation}
\label{proof.stability.2}
\| \partial_{t}(1/\sqrt{n_{\tau}})\|_{L^{\infty}(0,T;W^{-1,1}(\Omega))}\leq C,
\end{equation}
with $C$ a constant independent of $\eps$. Then \eqref{proof.stability.2} and \eqref{eq.proof.1} allow us to apply the Aubin-Simon's Lemma and to obtain the almost everywhere convergence of $1/\sqrt{n_\tau}$.

Finally, we have 
$$\nabla(n_{\tau}u_{\tau})=n_{\tau}\nabla u_{\tau}+ u_{\tau} \nabla n_{\tau}=\sqrt{n_{\tau}} \sqrt{n_{\tau}}\nabla u_{\tau}+2 \, \sqrt{n_{\tau}} \,  u_{\tau} \nabla \sqrt{n_{\tau}} \in L^2(0,T; L^1(\Omega)).$$
Therefore $n_{\tau}u_{\tau} \in L^2(0,T;W^{1,1}(\Omega))$. Using the momentum equation, we can get
information on $\partial_t(n_{\tau}u_{\tau})$ and therefore through Aubin-Lions-Simon's Lemma convergence almost everywhere of $n_{\tau}u_{\tau}$. 
With \eqref{eq.prop.2} and since $\sqrt{n_{\tau}} u_{\tau}$ converges almost everywhere then $(\sqrt{n_{\tau}} u_{\tau})^2$ converges strongly in $L^1(0,T;L^1(\Omega))$ and therefore  $\sqrt{n_{\tau}} u_{\tau}$ converges strongly in  $L^2(0,T;L^2(\Omega))$.
\cqfd

We can pass to the limit $\delta \rightarrow 0$ thanks to the strongly convergence on $\sqrt{n_\tau} u_\tau$ and the weakly convergence on $\nabla \sqrt{n_\tau}$ because
we write
$$n_\tau u_\tau=\sqrt{n_\tau}  \, \sqrt{n_\tau} u_\tau,$$
$$n_\tau u_\tau \otimes u_\tau=\sqrt{n_\tau} u_\tau \otimes \sqrt{n_\tau}  u_\tau,$$
	and
$$<n_\tau D(u_\tau),  \nabla \phi>=< D(n_\tau u_\tau),\nabla \phi>-<(u_{\tau})_j \partial_i n_\tau, \partial_i \phi_j>-<(u_\tau)_i \partial_j n_\tau,\partial_i \phi_j>$$
$$=-<\sqrt{n_\tau}\sqrt{n_\tau} u_\tau, D \nabla \phi>-<\sqrt{n_\tau} (u_n)_j \partial_i \sqrt{n_\tau},\partial_i \phi_j>-<\sqrt{n_\tau} (u_\tau)_i \partial_j \sqrt{n_\tau}, \partial_i \phi_j>.$$


\section{Low Planck limit}\label{sec.limit}

\vskip0.2cm
Let us consider $(n^\eps,u^\eps)_\eps$ a sequence of solutions of \eqref{eq2} in the sense of Definition \ref{def.solfaible1}. The goal of this section is to prove Theorem \ref{thm.convergence} {\it i.e.} that, up to a subsequence, $(n^\eps,u^\eps)_\eps$ tends to $(n^0,u^0)$ solution of \eqref{eql} in the sense of Definition \ref{def.solfaible2} when $\eps$ tends to 0. 

Using a priori estimates showed in Section 3 we can prove the same convergence as in the stability Section 4.2. The unique difference concerns the strong convergence of $\sqrt{n^\eps}$ which is here only in $L^2(0,T;L^2(\Omega))$ instead of $L^2(0,T;H^1(\Omega))$. Indeed this last one was obtain using an estimate which badly depend on $\eps$. 

\begin{proposition}\label{prop.lim.eps}
Under the hypothesis of Theorem \ref{thm.convergence}, we have,  when $\eps$ tends to $0$, 
\begin{align*}
&\sqrt{n^{\eps}} \to \sqrt{n}, \quad \hbox{strongly in} \ L^{2}(0,T;L^{2}(\Omega)),\\
&p(n^{\eps}) \to p(n), \quad \hbox{strongly in} \ L^{1}(0,T;L^{1}(\Omega)), \\
&p_{c}(n^{\eps}) \to p_{c}(n), \quad \hbox{strongly in} \ L^{1}(0,T;L^{1}(\Omega)),\\
&1/\sqrt{n^{\eps}} \to 1/\sqrt{n}, \quad \hbox{almost everywhere}, \\
&\sqrt{n^{\eps}}\,u^{\eps} \to \sqrt{n}\,u, \quad \hbox{strongly in} \ L^{2}(0,T;L^{2}(\Omega)).
\end{align*}
\end{proposition}

\pr
To prove the first convergence we use the estimate on $\nabla \sqrt{n}$ in $L^2(0,T;L^2(\Omega))$, the estimate on $\partial_t \sqrt{n}$ in $L^2(0,T;H^{-1}(\Omega))$ (obtained in Section 4.2) and the Aubin-Simon's  Lemma with $H^1(\Omega) \subset\subset L^2(\Omega) \subset H^{-1}(\Omega)$. We refer the reader to Section 4.2 for the esta\-blishment of the others convergences since they are completely similar.
\cqfd

Convergences of Proposition \ref{prop.lim.eps} allow us to pass to the limit $\eps$ tends to 0 in the second and third integrals of the left hand side and the last one of the right hand-side of the weak formulation \eqref{eq.solfaible1}. Using the fact that $\sqrt{n}$ lies in $L^2(0,T;L^2(\Omega))$ and $\nabla \sqrt{n}$ in $L^\infty(0,T;L^2(\Omega))$ (due to \eqref{eq.lemme2.2}) we can show that:
$$2 \, \eps^2  \,   \int_0^T \int_\Omega \sqrt n \,  \nabla \sqrt n \cdot \nabla {\rm div} \phi 
  +  4 \eps^2  \int_0^T \int_\Omega   \nabla\sqrt n \otimes \nabla\sqrt n : \nabla \phi \leq C\eps^2,$$
with $C$ a constant independent of $\eps$. Then these two integrals go to $0$ when $\eps$ tends to $0$, and we obtain the result.

\section*{Acknowledgements}
The authors would like to thank Prof. Didier Bresch.  All discussions with him were very helpful to achieve this paper.

\end{document}